\documentclass[11pt]{amsart}
\title{Selective screenability in topological groups
\footnote{\lowercase{{\bf {\uppercase{K}}ey words and phrases:} {\uppercase{H}}aver property, selective screenability, {\uppercase{H}}urewicz property, finitary {\uppercase{H}}aver property, countable dimensional, selection principle, c-groupable cover.\\
{{\bf{\uppercase{S}}ubject {\uppercase{C}}lassification:} {\uppercase{P}}rimary 54{\uppercase{D}}20, 54{\uppercase{D}}45, 55{\uppercase{M}}10; {\uppercase{S}}econdary 03{\uppercase{E}}20.}
}}
}

\author{Liljana Babinkostova}
\date{}
\usepackage{amssymb}
\usepackage{amsfonts}
\newcommand{\op}{\mathcal{O}}
\newcommand{\Onbd}{{\mathcal O}_{\sf nbd}}
\newcommand{\Omnbd}{{\Omega}_{\sf nbd}}
\newcommand{\sone}{{\sf S}_1}
\newcommand{\sfin}{{\sf S}_{fin}}
\newcommand{\Sc}{{\sf S}_c}
\newcommand{\Gc}{{\sf G}_c}
\newcommand{\naturals}{{\mathbb N}}
\newcommand{\reals}{{\mathbb R}}

\newcommand{\epf}{\diamondsuit}
  
\newtheorem{theorem}{Theorem}
\newtheorem{lemma}[theorem]{Lemma}
	
\newtheorem{corollary}[theorem]{Corollary}	
\newtheorem{problem}{Problem}
\newtheorem{conjecture}[problem]{Conjecture}
\begin{document}
\begin{abstract} 
We examine the selective screenability property in topological groups. In the metrizable case we also give characterizations of  $\Sc(\Onbd,\op)$ and {\sf Smirnov}-$\Sc(\Onbd,\op)$  in terms of the Haver property and finitary Haver property respectively relative to left-invariant metrics. We prove theorems stating conditions under which $\Sc(\Onbd,\op)$ is preserved by products.  Among metrizable groups we characterize the countable dimensional ones by a natural game.
\end{abstract}
\maketitle
\section{Definitions and notation}

Let $G$ be topological space. We shall use the notations:
\begin{itemize}
\item{$\mathcal{O}$:} {The collection of open covers of $G$.}
\end{itemize}
An open cover $\mathcal{U}$ of a topological space $G$ is said to be
\begin{itemize}
\item {an $\omega$-cover if $G$ is not a member of $\mathcal{U}$, but for each finite subset $F$ of $G$ there is a $U\in \mathcal{U}$ such that $F\subset U$. The symbol $\Omega$ denotes the collection of $\omega$ covers of $G$.}
\item { \emph{groupable} if there is a partition $\mathcal{U} =\cup_{n< \infty}\mathcal{U}_n$, where each $\mathcal{U}_n$ is finite, and for each $x\in G$ the set $\{n:x\not\in\cup\mathcal{U}_n\}$ is finite. The symbol $\mathcal{O}^{gp}$ denotes the collection of groupable open covers of the space.}
\item{\emph{large} if each element of the space is contained in infinitely many elements of the cover. The symbol $\Lambda$ denotes the collection of large covers of the space.}
\item{\emph{c-groupable} if there is a partition $\mathcal{U} =\cup_{n< \infty}\mathcal{U}_n$, where each $\mathcal{U}_n$ is pairwise disjoint and each $x$ is in all but finitely many $\cup{\mathcal{U}_n}$. The symbol $\mathcal{O}^{cgp}$ denotes the collection of c-groupable open covers of the space.}
\end{itemize}

Now let $(G,*)$ be a topological group with identity element $e$. We will assume that $G$ is not compact. For $A$ and $B$ subsets of $G$, $A*B$ denotes $\{a*b:a\in A,\,b\in B\}$. We use the notation $A^2$ to denote $A*A$, and for $n>1$, $A^n$ denotes $A^{n-1}*A$. For a neighborhood $U$ of $e$, and for a finite subset $F$ of $G$ the set $F*U$ is a neighborhood of the finite set $F$. Thus, the set $\{F*U:F\subset G \mbox{ finite}\}$ is an $\omega$-cover of $G$, which is denoted by the symbol $\Omega(U)$. The set
\[
  \Omega_{nbd} = \{\Omega(U):U \mbox{ a neighborhood of }e\}
\]  
is the set of all such $\omega$-covers of $G$. 

The set $\mathcal{O}(U) = \{x*U:x\in G\}$ is an open cover of $G$. The symbol
\[
  \mathcal{O}_{nbd} = \{\mathcal{O}(U): U \mbox{ a neighborhood of } e\}
\]
denotes the collection of all such open covers of $G$. Selection principles using these open covers of topological groups have been considered in several papers, including \cite{lbgps}, \cite{lblkms}, \cite{MT} and \cite{BT}, where information relevant to our topic can be found.
Now we describe the relevant selection principles for this paper. Let $S$ be an infinite set, and let $\mathcal{A}$ and $\mathcal{B}$ be collections of families of subsets of $S$.\\

The selection principle $\Sc(\mathcal{A},\mathcal{B})$, introduced in \cite{lbsc}, is defined as follows:
\begin{quote}
For each sequence $(A_n:n<\infty)$ of elements of the family $\mathcal{A}$ there exists a sequence $(B_n:n<\infty)$ such that for each $n$ $B_n$ is a pairwise disjoint family refining $A_n$, and $\bigcup_{n<\infty}B_n$ is a member of the family $\mathcal{B}$.
\end{quote}
The selection principle ${\sf Smirnov}-\Sc(\mathcal{A},\mathcal{B})$ is defined as follows:
\begin{quote}
For each sequence $(A_n:n<\infty)$ of elements of the family $\mathcal{A}$ there exists a positive integer $k<\infty$ and a sequence $(B_n:n\leq k)$  where each $B_n$ is a pairwise disjoint family of open sets refining $A_n$, $n\leq k$ and $\bigcup_{j\leq k}B_j$ is a member of the family $\mathcal{B}$.
\end{quote}

The metrizable space $X$ is said to be \emph{Haver} \cite{h} with respect to a metric $d$ if there is for each sequence $(\epsilon_n:n<\infty)$ of positive reals a sequence $(\mathcal{V}_n:n<\infty)$ where each $\mathcal{V}_n$ is a pairwise disjoint family of open sets, each of $d$-diameter less than $\epsilon_n$, such that $\bigcup_{n<\infty}\mathcal{V}_n$ is a cover of $X$. And it is said to be \emph{finitary Haver} \cite {borst} with respect to the metric $d$ if there is for each sequence $(\epsilon_n:n<\infty)$ a positive integer $k$ and a sequence $(\mathcal{V}_n:n\leq k)$ where each $\mathcal{V}_n$ is a pairwise disjoint family of open sets, each of diameter less than $\epsilon_n$, such that $\bigcup_{n\le k}\mathcal{V}_n$ is a cover of $X$.\\

\section{Selective screenability and $\Sc(\Onbd,\op)$}

Recent investigations into the Haver property and its relation to the selective screenability property $\Sc(\op,\op)$ revealed that the Haver property is weaker than selective screenability. E. and R. Pol has reported the following nice characterizations of $\Sc(\op,\op)$ in terms of the Haver property:
\begin{theorem}[\cite {erpol2}]\label{Pol} 
Let (X,d) be a metrizable space. The following are equivalent:
\begin{enumerate}
\item{X has property $\Sc(\op,\op)$.}
\item{X has the Haver property in all equivalent metrics.}
\end{enumerate}
\end{theorem}

For a topological group the property $\Sc(\op,\op)$ is stronger than $\Sc(\Onbd,\op)$. This is in part seen by comparing Theorem \ref{Pol} with the following result:  
\begin{theorem}\label{scandhavergroups} Let $(G,*)$ be a metrizable group. The following are equivalent:
\begin{enumerate}
\item{The group has property $\Sc(\op_{{\sf nbd}},\op)$.}
\item{The group has the Haver property in all equivalent left invariant metrics.}
\end{enumerate}
\end{theorem}

In the proof of Theorem \ref{scandhavergroups} we use the following result of Kakutani:
\begin{theorem}[\cite{HR}]\label{K} Let $(U_k:k<\infty)$ be a sequence of subsets of the topological group $(H,*)$ where $\{U_k:k<\infty\}$ is a neighborhood basis of the identity element $e$ and each $U_k$ is symmetric \footnote{$U_k$ is symmetric if $U_k = U_k^{-1}$},
and for each $k$ also $U^2_{k+1}\subseteq U_k$. Then there is a left-invariant metric $d$ on $H$ such that
\begin{enumerate}
\item{$d$ is uniformly continuous in the left uniform structure on $H\times H$.}
\item{If $y^{-1}*x\in U_k$ then $d(x,y)\leq (\frac{1}{2})^{k-2}$.}
\item{If $d(x,y) < (\frac{1}{2})^k$ then $y^{-1}*x\in U_k$.}
\end{enumerate}
\end{theorem}
In the above theorem $V^2$ denotes $\{a*b:\, a,\, b\in V\}$. For $n>1$ a positive integer the symbol $V^n$ is defined similarly. 

And now the proof of Theorem \ref{scandhavergroups}:\\
{\bf Proof:} $1\Rightarrow 2$: Let $d$ be a left-invariant metric of $G$ and let $(\epsilon_n:n<\infty)$ be a sequence of positive real numbers. For each $n$ choose an open neighborhood  $U_n$ of  the identity element $e$ of $G$ with $diam_d(U_n) <\epsilon_n$ and put  $\mathcal{U}_n=\mathcal{O}(U_n)$. Then $\{\mathcal{U}_n:n<\infty\}$ is a sequence from $\mathcal{O}_{\sf nbd}(U)$. Apply $\Sc(\op_{{\sf nbd}},\op)$. For each $n$ there is a pairwise disjoint family $\mathcal{V}_n$ of open sets refining $\mathcal{U}_n$ such that $\bigcup_{n<\infty}\mathcal{V}_n$ is an element of $\mathcal{O}$. Now for each $n$, for $V\in\mathcal{V_n}$ there is an $x\in G$ with $V\subseteq {x}* U_n$. But then $diam_d(V)\leq diam_d({x}* U_n)=diam_d(U_n)\leq\epsilon_n$. Thus the $\mathcal{V}_n$'s witnesses Haver's property for the given sequence of $\epsilon_n$'s.\\
$2\Rightarrow 1$: Let $\mathcal{U}_n=\mathcal{O}(U_n)$, $n<\infty$ be given. For each $n$ choose a neighborhood $V_n$ of the identity element $e$ in $G$ such that:  
\begin{enumerate}
\item{For all $n$, $V_n\subset U_n$.}
\item{For all $n$, $V_n*V_n\subset V_{n-1}$.}
\item{$\{V_n:n<\infty\}$is a neighborhood basis of the identity $e$.}
\end{enumerate}
By Kakutani's theorem choose a left invariant metric $d$ so that for every $n$:
\begin{enumerate}
\item{If $y^{-1}*x\in V_n$ then $d(x,y)\leq (\frac{1}{2})^{n-2}$.}
\item{If $d(x,y) < (\frac{1}{2})^n$ then $y^{-1}*x\in V_n$.}
\end{enumerate}

For each $n$, put $\epsilon_n=(\frac{1}{2})^{n}$. Since $G$ has the Haver property with respect to $d$, choose for each $n$ a pairwice disjoint family $\mathcal {V}_n$ of open sets such that: 
\begin{enumerate}
\item{For each $n$ and for each $V\in\mathcal{V}_n$, $diam_d(V)<\epsilon_n$.}
\item{$\bigcup_{n<\infty}\mathcal{V}_n$ covers $G$.}
\end{enumerate}
Then for every $n$ and for every $V\in\mathcal{V}_n$, there is and $x_V$ with $V\subseteq x*V_n\subseteq x_V*U_n\in\mathcal{U}_n$ and so $\mathcal{V}_n$ refines $\mathcal{U}_n$. But then $\mathcal{V}_n$ witness $\Sc(\op_{{\sf nbd}},\op)$ for $\{\mathcal{U}_n:n<\infty\}$. $\epf$

Using the similar ideas one can prove the following:
\begin{theorem}\label{smirscandfinhavergroups} Let $(G,*)$ be a metrizable group. The following are equivalent:
\begin{enumerate}
\item{The group has property ${\sf Smirnov}-\Sc(\op_{{\sf nbd}},\op)$.}
\item{The group has the finitary Haver property in all equivalent left invariant metrics.}
\end{enumerate}
\end{theorem}

One may ask when the properties $\Sc(\op,\op)$ and $\Sc(\Onbd,\op)$ are equivalent in a topological group. We do not have a complete answer. The Hurewicz property gives a condition: A topological space $G$ has the \emph{Hurewicz property} if for each sequence $\mathcal{U}_n, n<\infty$ of open covers of $X$ there is a sequence $\mathcal{F}_n, n<\infty$ of finite sets such that each $\mathcal{F}_n\subset\mathcal{U}_n$, and for each $x\in G$, the set $\{n:x\not\in\cup\mathcal{F}_n\}$ is finite.  

\begin{theorem}\label{hurewiczgroups}
Let $(G,*)$ be a topological group with the Hurewicz property. Then $\Sc(\op_{{\sf nbd}},\op)$ is equivalent to $\Sc(\op,\op)$. 
\end{theorem}
{\flushleft{{\bf Proof:}}} Let $(G,*)$ be a topological group. It is clear that $\Sc(\op,\op)$ implies $\Sc(\Onbd,\op)$. For the converse implication, assume the group has property $\Sc(\Onbd,\op)$. Let $(\mathcal{U}_n:n<\infty)$ be a sequence of open covers of $G$. For each $n$, and each $x\in G$ choose a neighborhood $V(x,n)$ of the identitity $e$ such that $x*V(x,n)^4$ is a subset of some $U$ in $\mathcal{U}_n$. Put $\mathcal{H}_n=\{x*V(x,n):x\in G\}$.
Apply the Hurewicz property to the sequence $(\mathcal{H}_n:n< \infty)$. For each $n$ choose a finite $\mathcal{F}_n\subset \mathcal{H}_n$ such that for each $g\in G$, the set $\{n:g\not\in\cup\mathcal{F}_n\}$ is finite. Write $\mathcal{F}_n=\{{x_i}^n*V({x_i}^n,n):i\in I_n\}$ and $I_n$ is finite. For each $n$, define $V_n=\bigcap_{i\in I_n}V({x_i}^n,n)$ a neighborhood of the identity $e$. Choose a partition $\naturals= \bigcup_{k<\infty}J_k$ where each $J_k$ is infinite, and for $l\neq k$, $J_l\cap J_k = \emptyset$. For each $k$, apply $\Sc(\Onbd,\op)$ to the sequence $(\mathcal{O}(V_n):n\in J_k)$. For each $n\in J_k$ find a pairwise disjoint family $\mathcal{S}_n$ of open sets such that $\mathcal{S}_n$ refines $\mathcal{O}(V_n)$ and $\bigcup_{n\in J_k} \mathcal{S}_n$ covers $G$. For each $n$ define $\mathcal{V}_n=\{S\in \mathcal{S}_n:(\exists U\in \mathcal{U}_n)(S\subseteq U)\}$. Since $\mathcal{V}_n\subset \mathcal{S}_n$, $\mathcal{V}_n$ is pairwise disjoint and refines $\mathcal{U}_n$. We will show that $\bigcup_{n<\infty}\mathcal{V}_n$ covers $G$. Pick any $g\in G$. Fix $N_g$ so that for all $n\ge N_g$, $g\in\cup\mathcal{F}_n$. Pick $k_g$ so large that $min(J_{k_g})\ge N_g$. Pick $m\in J_{k_g}$ with $g\in\cup\mathcal{S}_m$. Pick $J\in\mathcal{S}_m$ with $g\in J$. We will show that $J\in\mathcal{V}_m$. We have that $g\in\cup\mathcal{F}_m$, so pick $i\in I_m$ with $g\in {x_i}^m*V({x_i}^m,m)$. Since $J\in\mathcal{S}_m$, also pick $h_m$ so that $J\subseteq h_m * V_m = h_m\times(\bigcap_{i\in I_m}V({x_i}^m,m)\subseteq h_m * V({x_i}^m,m)$. We have that $g={x_i}^m*z_g=h_m*t_g$ for some $z_g,t_g \in V({x_i}^m,m)$. So $h_m={x_i}^m*z_g*{t_g}^{-1}$. Now consider any $y\in J$. Choose $t_y\in V({x_i}^m,m)$ with $y=h_m*t_y$. But then $y= {x_i}^m*(z_g*{t_g}^{-1}*t_y*e)\in {x_i}^m*V^4({x_i}^m,m)\subseteq U$, for some $U\in\mathcal{U}_m$. So we have that $J\in\mathcal{V}_m$ and $g\in J$.$\epf$

The symbol $\sone(\mathcal{A},\mathcal{B})$ denotes the statement that there is for each sequence $(O_n:n<\infty)$ of elements of $\mathcal{A}$ a sequence $(T_n:n<\infty)$ such that for each  $n$ $T_n\in O_n$, and $\{T_n:n<\infty\}\in\mathcal{B}$. 
A topological group $(G,*)$ is said to be a \emph{Hurewicz-bounded} group if it satisfies the selection principle
$\sone(\Omega_{nbd},\mathcal{O}^{gp})$. In \cite{lbsc} was shown that $\Sc(\mathcal{O},\mathcal{O})$ is equivalent to $\Sc(\Omega,\mathcal{O})$. The analogous equivalence doesn't hold in topological groups: 
\begin{theorem} 
$\Sc(\Omnbd,\op)$ does not imply $\Sc(\mathcal{O}_{\sf nbd},\mathcal{O})$.
\end{theorem}
{\bf Proof:} Let $(C,*)$ be the unit circle in the complex plane with complex multiplication. It is a compact metrizable group embedding the unit interval $[0,1]$ as a subspace. Since $(C^{\naturals},*)$ is a compact group it has the Hurewicz property, so is Hurewicz bounded.  Also $\reals$, the real line with addition, is a Hurewicz-bounded topological group. Thus the product group $\reals\times C^{\naturals}$ is Hurewicz bounded, so has the property $\sone(\Omnbd,\op)$, and so has $\Sc(\Omnbd,\op)$. But $[0,1]^{\naturals}$ embeds as closed subspace into $\reals\times C^{\naturals}$, and $[0,1]^{\naturals}$ does not have the property $\Sc(\op,\op)$. Thus the topological group $\reals\times C^{\naturals}$ does not have $\Sc(\op,\op)$, and as it has the Hurewicz property, Theorem \ref{hurewiczgroups} implies it is not $\Sc(\Onbd,\op)$.$\epf$ 

The symbol $\sfin(\mathcal{A},\mathcal{B})$ denotes the statement that there is for each sequence $(O_n:n<\infty)$ of elements of $\mathcal{A}$ a sequence $(T_n:n<\infty)$ of finite sets such that for each  $n$ $T_n\subseteq O_n$, and $\cup\{T_n:n<\infty\}\in\mathcal{B}$. It was shown in \cite{KS} that $\sfin(\Omega,\op^{gp})$ is equivalent to the Hurewicz property. And it is well known that $\sfin(\op,\op)$ is the Menger property, which is equivalent to $\sfin(\Omega,\op)$. A topological group is said to be a \emph{Menger bounded} group if it has the property $\sone(\Omnbd,\op)$. 

By how much can the requirement that $(G,*)$ has the Hurewicz property be weakened in Theorem \ref{hurewiczgroups}? Natural possibilities include the Menger property, Menger boundedness or Hurewicz boundedness. In light of interesting recent examples of E. and R. Pol - \cite{erpol}, \cite{erpol2} we conjecture that none of these weakenings is enough:
\begin{conjecture}\label{mengerbounded} There is a metrizable Menger bounded topological group which has the property $\Sc(\Onbd,\op)$, but not the property $\Sc(\op,\op)$.
\end{conjecture}
\begin{conjecture}\label{hurewiczbounded} There is a metrizable Hurewicz bounded topological group which has the property $\Sc(\Onbd,\op)$, but not the property $\Sc(\op,\op)$.
\end{conjecture}
\begin{conjecture}\label{menger} There is a metrizable topological group which has the Menger property and property $\Sc(\Onbd,\op)$, but not the property $\Sc(\op,\op)$.
\end{conjecture}
It is clear that Conjecture \ref{menger}$\Rightarrow$ Conjecture \ref{mengerbounded} and Conjecture \ref{hurewiczbounded} $\Rightarrow$ Conjecture \ref{mengerbounded}. It may be that Conjecture \ref{hurewiczbounded} is independent of the Zermelo-Fraenkel axioms. Recently  E. and R. Pol showed that CH implies Conjecture \ref{menger}.

\section{Products}

E. Pol showed in \cite{epol86} that there exist a zerodimensional subset $Y$ of the real line and a separable metric space $X$ and  such that $X$ has the property $\Sc(\op,\op)$ in all finite powers, but $X\times Y$ does not have $\Sc(\op,\op)$. This failure does not happen for the group analogue:

\begin{theorem}
Let $(G,*)$ be a group satisfying $\Sc(\Onbd,\op)$. If $(H,*)$ is a group with property $\Sc(\Onbd,\mathcal{O}^{cgp})$, then $(G\times H,*)$ also has $\Sc(\Onbd,\op)$.
\end{theorem}
{\bf Proof:} For each $n$ let $\mathcal{U}_n$ be an element of $\Onbd(G\times H)$. Each $\mathcal{U}_n$ is of the form $\mathcal{U}_n=\mathcal{O}(U_n)$ where $U_n$ is a neighborhood of the identity $(e_G,e_H)$ of $G\times H$. Pick $V_n\subset G$ a neighborhood of $e_G$, and $W_n\subset H$ a neighborhood of $e_H$ so that $V_n\times W_n\subseteq U_n$. Then $\mathcal{W}_n=\mathcal{O}(V_n\times W_n)$ is a refinement of $\mathcal{U}_n$, for all $n$.
Let $\mathcal{H}_n=\mathcal{O}(W_n)\in\mathcal{O}_{\sf nbd}$. Apply $\Sc(\Onbd,\mathcal{O}^{cgp})$ to the sequence $(\mathcal {H}_n:n<\infty)$. For each $n$ find a finite pairwise disjoint refinement $\mathcal{K}_n$ of $\mathcal{H}_n$ so that each $x$ is in all but finitely many of $\bigcup \mathcal{K}_n$. Next, for each $n$ put $\mathcal{G}_n=\mathcal{O}(V_n)\in \mathcal{O}_{nbd}$. Apply $\Sc(\Onbd,\mathcal{O})$ to the sequence $(\mathcal {G}_n:n< \infty)$. For each $n$ choose pairwise disjoint $\mathcal{J}_n$ that refines $\mathcal{G}_n$ so that $\bigcup \mathcal{J}_n$ is a large open cover of $G$.
For each $n$ define $\mathcal{V}_n=\{J\times K: J\in \mathcal{J}_n, K\in\mathcal{K}_n\}$.\\
{\bf Claim 1:} $\mathcal{V}_n$ refines $\mathcal{W}_n$: For $J\in \mathcal{J}_n$ and $K\in\mathcal{K}_n$ there is an element $g\in G$ and $h\in H$ such that $J\subseteq g* V_n$ and $K\subseteq g*W_n$. But then $J\times K\subseteq g* V_n \times h*W_n\in \mathcal{W}_n$. \\
{\bf Claim 2:} $\mathcal {V}_n$ is pairwise disjoint: Let $J_1\times K_1$ and $J_2\times K_2$ be elements of $\mathcal {V}_n$ with $J_1\times K_1 \neq J_2\times K_2$. If $J_1\neq J_2$ then $J_1\cap J_2=\emptyset$ because the $\mathcal {J}_n$ is disjoint. So $(J_1\times K_1)\bigcap (J_2\times K_2)=\emptyset$.
Similarly, $(J_1\times K_1)\bigcap (J_2\times K_2)=\emptyset$ if $K_1\neq K_2$.\\
{\bf Claim 3:} $\bigcup\mathcal{V}_n$ covers $G\times H$. Consider $(g,h)$ as an element of $C\times H$. Since $\bigcup \mathcal{J}_n$ is a large cover of $G$ the set $S_1=\{ n:(\exists J\in\mathcal{J}_n)(g\in J)\}$ is infinite and there is an $N$ such that $S_2=\{ n:(\exists K\in\mathcal{K}_n)(h\in K)\}\supseteq \{n: n\ge N\}$.
Pick an $n\in S_1\cap S_2$. Pick $J\in\mathcal{J}_n$ with $g\in J$ and $K\in\mathcal{K}_n$ with $h\in K$. Then $(g,h)\in J\times K\in\mathcal{V}_n$.$\epf$

\begin{corollary}\label{scgroupproducts} Let $(G,*_1)$ and $(H,*_2)$ be metrizable topological groups such that $(G,*_1)$ has $\Sc(\Onbd,\op)$ and  $H$ is zero-dimensional. Then $(G\times H,*)$ is a group with property $\Sc(\op_{{\sf nbd}},\op)$.
\end{corollary}
{\bf Proof:} We show that $(H,*)$ has $\Sc(\Onbd,\op^{cgp})$. The reason for this is that since $H$ is zerodimensional, each open cover of it has a refinement by a disjoint open cover. Thus for a given sequence $(\mathcal{U}_n:n<\infty)$ from $\op_{{\sf nbd}}$ for $H$ we can choose for each $n$ a disjoint open refinement $\mathcal{V}_n$ which covers $H$. Clearly $\cup_{n<\infty}\mathcal{V}_n$ is $c$-groupable.$\epf$

To illustrate: Let ${\mathbb P}$ denote the set of irrational numbers. E. Pol has shown under CH\footnote{For a new proof using a weaker hypothesis, see \cite{erpol} and \cite{erpol2}.} that there is a metrizable space $X$ with property $\Sc(\op,\op)$ such that $X\times{\mathbb P}$ does not have $\Sc(\op,\op)$. Now ${\mathbb P}$ is homeomorphic to a closed subset of the zerodimensional group $({\mathbb Z}^{\naturals},+)$. Thus $X\times {\mathbb Z}^{\naturals}$ also does not have $\Sc(\op,\op)$. But for any topological group $(G,*)$ with property $\Sc(\Onbd,\op)$, the group $G\times {\mathbb Z}^{\naturals}$ still has $\Sc(\Onbd,\op)$.

Hattori, Yamada and independently Rohm, have proven the following product theorem for $\Sc(\op,\op)$:
\begin{theorem}[Hattori-Yamada, Rohm]\label{hyr} If $X$ is $\sigma$-compact and if $X$ and $Y$ both have the property $\Sc(\op,\op)$, then $X\times Y$ has the property $\Sc(\op,\op)$.
\end{theorem}
We shall prove an analogous theorem, Theorem \ref{hyrforgroups}, for topological groups. Since $\Sc(\Onbd,\op)$ is weaker than $\Sc(\op,\op)$ (see the remarks following Conjecture \ref{menger}), we are able to use a weaker restriction than $\sigma$-compact. We use the following result in our proof:

\begin{lemma}[\cite{lbms}]\label{scandhurewicz} 
The following statements are equivalent:
\begin{enumerate}
\item{$X$ has the Hurewicz property and property $\Sc(\mathcal{O},\mathcal{O})$.} 
\item{For each sequence $(\mathcal{U}_n:n<\infty)$ of open covers of $X$ there is a sequence $(\mathcal{V}_n:n<\infty)$ such that:
      \begin{enumerate}
      \item{Each $\mathcal{V}_n$ is a finite collection of open sets;}
      \item{Each $\mathcal{V}_n$ is pairwise disjoint;}
      \item{Each $\mathcal{V}_n$ refines $\mathcal{U}_n$;}
      \item{there is a sequence $n_1 < n_2 < \cdots < n_k < \cdots$ of positive integers such that each element of $X$ is in all but finitely many of the sets $\cup(\cup_{n_k\le j< n_{k+1}}\mathcal{V}_j)$.}
      \end{enumerate}
}
\end{enumerate}
\end{lemma}

\begin{theorem}\label{hyrforgroups} Let $(G,*)$ be a group which has property $\Sc(\mathcal{O}_{nbd},\mathcal{O})$ as well as the Hurewicz property. Then for any  topological group $(H,*)$ satisfying $\Sc(\mathcal{O}_{nbd},\Lambda)$, $G\times H$ also satisfies $\Sc(\mathcal{O}_{nbd},\mathcal{O})$.
\end{theorem} 
{\bf Proof:}  Let $(\mathcal{O}(U_n\times V_n):n<\infty)$ be a sequence of $\mathcal{O}_{nbd}$-covers of $G\times H$. Then each $\mathcal{O}(U_n)$ is an $\mathcal {O}_{nbd}$-cover of $G$ and each $\mathcal{O}(V_n)$ is an $\mathcal {O}_{nbd}$-cover of $H$. 

Since $(G,*)$ has the Hurewicz property and $\Sc(\Onbd,\op)$, it has by Theorem \ref{hurewiczgroups} the property $\Sc(\op,\op)$. 
Letting $(\Onbd(U_n):n<\infty)$ be the sequence of open covers in $(2)$ of Lemma \ref{scandhurewicz}, let $(\mathcal{V}_n:n<\infty)$ be the corresponding sequence provided by $(2)$ of that lemma, and fix $n_1 < n_2 < \cdots < n_{k+1} < \cdots$ as there.

For each $k$ define $W_k = \cap_{n_k\le j<n_{k+1}}V_j$. Then consider the sequence $(\Onbd(W_k):k<\infty)$ for $H$. 
Since $(H,*)$ has property  $\Sc(\Onbd,\Lambda)$ choose for each $k$ a pairwise disjoint refinement $\mathcal{R}_k$ of $\Onbd(W_k)$, consisting of open sets, such that each $h\in H$ is contained in infinitely many of the sets $\cup\mathcal{R}_k$. Notice that for each $k$, $\mathcal{R}_k$ is a disjoint refinement of each $\Onbd(V_j)$ for $n_k\le j<n_{k+1}$.

For each $j$ define $\mathcal{K}_j$ as follows: Find $k$ with $n_k\le j< n_{k+1}$ and put 
\[
  \mathcal{K}_j = \{V\times R: V\in\mathcal{V}_j \mbox{ and }R\in\mathcal{R}_k\}.  
\]

{\flushleft{\bf Claim 1:}} $\mathcal{K}_j$ is a refinement of $\Onbd(U_j\times V_j)$:\\
\underline{Proof:} Consider $V\times R\in\mathcal{K}_j$: Since $V\in\mathcal{V}_j$, choose a member $A_j$ of $\Onbd(U_j)$ with $V\subset A_j$. Choose $g_j\in G$ with $A_j = g_j*U_j$. Next, since $R\in\mathcal{R}_k$, choose a $B_k\in\Onbd(W_k)$ with $R\subseteq B_k$. Choose $h_j\in H$ with $B_k = h_j*W_k$. Then in particular we have $R\subseteq B_k \subseteq h_j*V_j$. But this implies that $V\times R\subset (g_j,h_j)*(U_j\times V_j)$, an element of $\Onbd(U_j\times V_j)$.\\ 

{\flushleft{\bf Claim 2:}} $\mathcal{K}_j$ is a disjoint family of open sets:\\
\underline{Proof:} This is clear.\\

{\flushleft{\bf Claim 3:}} $\cup_{j<\infty}\mathcal{K}_j$ is a cover of $G\times H$:\\
\underline{Proof:} To see this, consider $(g,h)\in G\times H$. Choose $N$ so large that for each $k\ge N$ we have $g\in \cup(\cup_{n_k\le j < n_{k+1}}\mathcal{V}_j)$. Then choose a $k>N$ with $h\in \cup\mathcal{R}_k$. It follows that for a $j$ with $n_k\le j < n_{k+1}$ we have $(g,h)$ in $\cup\mathcal{K}_j$.\\
This completes the proof. $\diamondsuit$

\begin{theorem}\label{largecovers} Let $(G,*)$ be a metrizable topological group with no isolated points. Then $\Sc(\Onbd,\op)$ is equivalent to $\Sc(\mathcal{O}_{\sf nbd},\Lambda)$.
\end{theorem}
{\bf Proof:} Let $(\mathcal{O}(U_n):n<\infty)$ be a sequence in $\Onbd(G)$. Choose a sequence $\epsilon_n:n<\infty)$ such that $\epsilon_{i}>\epsilon_{i+1}$ for all $i<\infty$ and $diam_d(U_1\cap U_2\cap\cdots\cap U_n)>\epsilon_n$ for all $n$. Define $(\mathcal{O}(V_n):n<\infty)$ such that $diam_d(V_i)=\epsilon_i$ for $i=1,2,\cdots,n$. Write $\naturals= \bigcup_{m<\infty}I_m$ where each $I_m$ is infinite, and for $m\neq k$, $I_m\cap I_k = \emptyset$. Apply $S_c(\Onbd,\op)$ to the sequence $(\mathcal{O}(V_n):n\in\ I_m)$ for all $m$. Let $T_n$ be a pairwise disjoint family refining $\mathcal{O}(V_n)$, $n\in I_m$  such that $\cup\{T_n:n\in I_m\}$ covers $G$. We will show that $\cup\{T_n:n\in\naturals\}$ is a large cover. Take an element $x\in G$ and pick $m_1\in I_1$ with $x\in\cup\ {T_{m_1}}$. Next, pick $W_1\in T_{m_1}$ with $x\in W_1$ and $N_1$ so large that for all $n\ge N_1$ we have $\epsilon_n < diam_d(W_1)$. Then pick  $i_2$ so large that the smallest element of $I_{i_2}$ is larger than $N_1$. Now choose $m_2\in I_{i_2}$ with $x\in \bigcup T_{m_2}$. Pick $W_2\in T_{m_2}$ with $x\in W_2$. Since $m_2\ge N_1$, $\epsilon _{m_2} < diam_d(W_1)$, and by definition of $\mathcal{O}(V_{m_2})$, $diam_d(W_2)\leq diam_d(V_{m_2})\leq \epsilon_{m_2} < diam_d(W_1)$. Next pick $N_2$ so large that for all $n\ge N_2$ we have $\epsilon_n < diam_d(W_2)$ and continue the same way as we did with $N_1$. Continuing like this we find $W_1, W_2, W_3, \cdots$ infinitely many distinct elements of $\cup\{T_n:n<\infty\}$ covering $x$.$\epf$

Note in particular that if for each $n$ $\mathcal{V}_n$ is a disjoint family of open sets, and if $\cup_{n<\infty}\mathcal{V}_n$ is a large cover of $G$, then for each $g\in G$ the set $\{n: g\in \cup\mathcal{V}_n\}$ is infinite. This is because for each $n$ there is at most one set in $\mathcal{V}_n$ that might contain $g$.

\begin{corollary}\label{hyrformgroups} Let $(G,*)$ be a group which has property $\Sc(\mathcal{O}_{nbd},\mathcal{O})$ as well as the Hurewicz property. Then for any metrizable topological group $(H,*)$ satisfying $\Sc(\mathcal{O}_{nbd},\op)$, $G\times H$ also satisfies $\Sc(\mathcal{O}_{nbd},\mathcal{O})$.
\end{corollary} 
{\bf Proof:} Use Theorems \ref{hyrforgroups} and \ref{largecovers}. $\epf$

\begin{corollary}\label{powersformgroups} Let $(G,*)$ be a metrizable group which has property $\Sc(\mathcal{O}_{nbd},\mathcal{O})$ as well as the Hurewicz property. Then all finite powers of $(G,*)$ have the property $\Sc(\mathcal{O}_{nbd},\mathcal{O})$.
\end{corollary} 
{\bf Proof:} Use Corollary \ref{hyrformgroups}. $\epf$

It is not clear that the full Hurewicz property is needed in Theorem \ref{hyrforgroups} or Corollaries \ref{hyrformgroups} and \ref{powersformgroups}: maybe Hurewicz-boundedness is enough.

\begin{problem} In Theorem \ref{hyrforgroups}, can we replace the requirement that $G$ has the Hurewicz property with the weaker requirement that $(G,*)$ has the property $\sone(\Omnbd,\mathcal{O}^{gp})$? 
\end{problem}

In light of  results of E. and R. Pol - \cite{erpol} - we conjecture that neither Menger boundedness, nor the Menger property is enough to obtain Theorem \ref{hyrforgroups}:
\begin{conjecture}\label{mengerbounded2} There is a metrizable Menger bounded group $(G,*)$ with property $\Sc(\Onbd,\op)$, such that $G^2$ is Menger bounded but does not have $\Sc(\Onbd,\op)$.
\end{conjecture}
\begin{conjecture}\label{mengergroup} There is a metrizable group $(G,*)$ which has the property $\Sc(\Onbd,\op)$, and $G^2$ has the Menger property but does not have $\Sc(\Onbd,\op)$.
\end{conjecture}

\section{Games}

The following game, denoted $\Gc(\mathcal{A},\mathcal{B})$, is naturally associated with $\Sc(\mathcal{A},\mathcal{B})$: Players ONE and TWO play as follows: They play an inning for each natural number $n$. In the $n$-th inning ONE first chooses $\mathcal{O}_n$, a member of $\mathcal{A}$, and then TWO responds with $\mathcal{T}_n$ refining $\mathcal{O}_n$. A play $(\mathcal{O}_1, \mathcal{T}_1, \cdots, \mathcal{O}_n, \mathcal{T}_n, \cdots)$ is won by TWO if $\cup_{n<\infty}\mathcal{T}_n$ is a member of $\mathcal{B}$; else, ONE wins. Versions of different length of this game can also be considered: For an ordinal number $\alpha$ let $\Gc^{\alpha}(\mathcal{A},\mathcal{B})$ be the game played as follows: in the $\beta$-th inning ($\beta<\alpha$) ONE first chooses $\mathcal{O}_{\beta}$, a member of $\mathcal{A}$, and then TWO responds with a pairwise disjoint $\mathcal{T}_{\beta}$ which refines $\mathcal{O}_{\beta}$. A play
\[
  \mathcal{O}_0,\mathcal{T}_0, \cdots, \mathcal{O}_{\beta}, \mathcal{T}_{\beta},\cdots \, \beta<\alpha
\]  
is won by TWO if $\cup_{\beta<\alpha}\mathcal{T}_{\beta}$ is a member of $\mathcal{B}$; else, ONE wins. Thus the game $\Gc(\mathcal{A},\mathcal{B})$ is $\Gc^{\omega}(\mathcal{A},\mathcal{B})$.

\begin{theorem} Let $(G,*)$ be a metrizable group. Then the following statements hold:
\begin{enumerate}
\item{If $dim(G)\leq n$ then TWO has a winning strategy in $\Gc^{n+1}(\mathcal{O}_{\sf nbd},\mathcal{O})$.}
\item{If TWO has a winning strategy in $\Gc^{n+1}(\mathcal{O}_{\sf nbd},\mathcal{O})$, then the $dim(G)\leq n$.}
\item{If $G$ is countable dimensional, then TWO has a winning strategy in $\Gc^{\omega}(\mathcal{O}_{\sf nbd},\mathcal{O})$.}
\item{If TWO has a winning strategy in $\Gc^{\omega}(\mathcal{O}_{\sf nbd},\mathcal{O})$, then $G$ is countable dimensional.}
\end{enumerate}
\end{theorem}
{\bf Proof:}We prove $3$ and $4$. The proofs of $1$ and $2$ are similar.\\
{\bf Proof of $3$:} 

Suppose that $G$ is countable dimensional. We define the following strategy for TWO: Write $G = \cup_{n<\infty}G_n$ where each $G_n$ is zero-dimensional. Let $\mathcal{U}$ be an element of $\mathcal{O}_{\sf nbd}$. For  $\mathcal{U}=\mathcal {O}(U)$ of $G$ and $n<\infty$, consider $\mathcal{U}$ as a cover of $G_n$. Since $G_n$ is zero-dimensional, find a pairwise disjoint family $\mathcal{V}_n$ of subsets of $G_n$ open in $G_n$ such that $\mathcal{V}_n$ covers $G_n$ and refines $\mathcal{O}(U)$. Choose a pairwise disjoint family $\sigma(\mathcal{U},n)$ refining $\mathcal{O}(U)$ such that each element $V$ of $\mathcal{V}_n$ is of the form $U\cap G_n$ for some $U\in\sigma(\mathcal{U},n)$. Now TWO plays as follows: In inning 1 ONE plays $\mathcal{U}_1$, and TWO responds with $\sigma(\mathcal{U}_1,1)$, thus covering $G_1$. When ONE has played $\mathcal{U}_2$ in the second inning TWO responds with $\sigma(\mathcal{U}_2,2)$, thus covering $G_2$, and so on. And in the $n$-th inning, when ONE has chosen the cover $\mathcal{U}_n$ of $G$ TWO responds with $\sigma(\mathcal{U}_n,n)$, covering $G_n$. This strategy evidently is a winning strategy for TWO.  \\
{\bf Proof of $4$:} Let $\sigma$ be a winning strategy for TWO. Choose a neighborhood basis $(U_n:n<\infty)$ of the identity element $e$ of $G$ so that $diam_d(U)<\frac{1}{n}$ for all $n$. Consider the plays of the game in which in each inning ONE chooses for some $n$ a cover $\mathcal{U}_n$ of $G$ of the form $\mathcal{O}(U_n)$.

Define a family $(C_{\tau}:\tau\in\, ^{<\omega}\naturals)$ of subsets of $G$ as follows:
\begin{enumerate}
\item{$C_{\emptyset} = \cap\{\cup \sigma(\mathcal{U}_n):n<\infty\}$;}
\item{For $\tau = (n_1,\cdots,n_k)$, $C_{\tau} = \cap\{\cup\sigma(\mathcal{U}_{n_1},\cdots,\mathcal{U}_{n_k},\mathcal{U}_n):n<\infty\}$}
\end{enumerate} 
Claim 1: $G = \cup\{C_{\tau}:\tau\in\, ^{<\omega}\naturals\}$:\\
For suppose on the contrary that $x\not\in \cup\{C_{\tau}:\tau\in\, ^{<\omega}\naturals\}$. Choose an $n_1$ such that $x\not\in\sigma(\mathcal{U}_{n_1})$. With $n_1,\cdots,n_k$ chosen such that $x\not\in\sigma(\mathcal{U}_{n_1},\cdots,\mathcal{U}_{n_k})$, choose an $n_{k+1}$ such that $x\not\in\sigma(\mathcal{U}_{n_1},\cdots,\mathcal{U}_{n_{k+1}})$, and so on. Then
\[
 \mathcal{U}_{n_1},\, \sigma(\mathcal{U}_{n_1}),\, \mathcal{U}_{n_2},\, \sigma(\mathcal{U}_{n_1},\mathcal{U}_{n_2}),\cdots
\]
is a $\sigma$-play lost by TWO, contradicting the fact that $\sigma$ is a winning strategy for TWO.\\
Claim 2: Each $C_{\tau}$ is zero-dimensional.\\
For consider an $x\in C_{\tau}$. Say $\tau = (n_1,\cdots,n_k)$. Thus, $x$ is a member of $\cap\{\cup\sigma(\mathcal{B}_{n_1},\cdots,\mathcal{U}_{n_k},\mathcal{U}_n):n<\infty\}$. For each $n$ choose a neighborhood $V_n(x)\in\sigma(\mathcal{U}_{n_1},\cdots,\mathcal{U}_{n_k},\mathcal{U}_n)$. Since for each $n$ we have $diam_d(V_n(x))<\frac{1}{n}$, the set $\{V_n(x)\cap C_{\tau}:n<\infty\}$ is a neighborhood basis for $x$ in $C_{\tau}$. Observe also that each $V_n(x)$ is closed in $C_{\tau}$ because: The set $V = \cup\sigma(\mathcal{U}_{n_1},\cdots,\mathcal{U}_{n_k},\mathcal{U}_n)\setminus V_n(x)$ is open in $G$ and so $C_{\tau}\setminus V_n(x) = C_{\tau}\cap V$ is open in $C_{\tau}$. Thus each element of $C_{\tau}$ has a basis consisting of clopen sets.
$\epf$

\section{Remarks and acknowledgment.}

Regarding Theorem \ref{scandhavergroups}: For a left invariant metric $d$ let $\mathcal{U}_d $ be a the family of sets $U_\epsilon$, $\epsilon > 0$ where we define $U_\epsilon =\{(x,y)\in G\times G: d(x,y)< \epsilon\}$.
The family $\mathcal{U}_d $ generates the left-uniformity of the topological group $G$. Refer to \cite{nb} Chapters III \S3 and IX \S3 and \cite {ek} Chapter 8.1 regarding these facts. Let $\mathcal{O}_d $ denote the collection of open covers of the form $\{U_\epsilon(x): x\in G\}$ where $U_\epsilon(x)= \{y: (x,y)\in U_\epsilon\}$. The referee pointed out that a third equivalence can be added in Theorem \ref{scandhavergroups}: 
\begin{quote}
{\tt For each left-invariant metric $d$, $S_c(\mathcal{O}_d,\op)$ holds.}
\end{quote}

Regarding Theorem \ref{hurewiczgroups}: There is a more general theorem. 
Let $\mathcal{U}$ be uniformity on $X$ generating the topology $\tau_{\mathcal{U}}$. For $V\subset X\times X$, define $V(x)=\{y\in X: (x,y)\in V\}$. We say that an open cover of $(X,\tau_\mathcal{U})$ is uniform with respect to $\mathcal{U}$ if it is of the form $\{V(x):x\in X\}$, for some $V\in \mathcal{U}$ . Define $\mathcal{O}_{\mathcal{U}}=\{\{V(x):x\in X\}:V\in\mathcal{U}\}$.
\begin{theorem}\label{uniformity} Let $\mathcal{U}$ be a uniformity generating the topology $\tau_{\mathcal{U}}$ on the set $G$. Assume that the topological space $(G,\tau_\mathcal{U})$ has the Hurewicz property. Then $S_c(\mathcal{O}_{\mathcal{U}},\op)$  is equivalent to $S_c(\op,\op)$.
\end{theorem}
The proof of this theorem is very similar to the proof of Theorem \ref{hurewiczgroups}. 

I thank the referee for the useful remarks and E. Pol and R. Pol for communicating their result on Conjecture \ref {menger} to me.

{\flushleft Contact Information:}\\
\begin{center}
\begin{tabular}{l}
Liljana Babinkostova      \\
Department of Mathematics \\
Boise State University    \\
Boise, ID 83725           \\
                          \\
e-mail:                    liljanab@diamond.boisestate.edu \\
fax:                       208 - 426 - 1356                \\
phone:                     208 - 426 - 2896
\end{tabular}
\end{center}
\end{document}